\documentclass[twosided,reqno]{amsart}

\usepackage{amsmath}
\usepackage{amssymb}
\usepackage{amsthm}

\theoremstyle{theorem}
\newtheorem{theorem}{\scshape Theorem }[section]
\newtheorem{lemma}[theorem]{\scshape Lemma}

\theoremstyle{definition}

\numberwithin{equation}{section}

\begin{document}

\title[Higher-order Bernoulli, Frobenius-Euler and Euler polynomials]{Higher-order Bernoulli, Frobenius-Euler and Euler polynomials}

\author{Taekyun Kim$^1$}
\address{$^1$ Department of Mathematics, Kwangwoon University, Seoul 139-701, Republic of Korea.}
\email{tkkim@kw.ac.kr}

\author{Dae San Kim$^2$}
\address{$^2$ Department of Mathematics, Sogang University, Seoul 121-742, Republic of Korea.}
\email{dskim@sogang.ac.kr}


\maketitle

\begin{abstract}
In this paper, we give some interesting identities of higher-order Bernoulli, Frobenius-Euler and Euler polynomials arising from umbral calculus. From our method of this paper, we can derive many interesting identities of special polynomials.
\end{abstract}

\section{Introduction}

For $\alpha \in {\mathbb{R}}$, the {\it{Bernoulli polynomials}} of order $\alpha$ are defined by the generating function to be
\begin{equation}\label{1}
\left(\frac{t}{e^t-1}\right)^{\alpha} e^{xt}=\sum_{n=0} ^{\infty}B_n ^{(\alpha)}(x)\frac{t^n}{n!},{\text{ (see [1,7,8,12,17,18])}}.
\end{equation}
In the special case, $x=0$, $B_n ^{(\alpha)} (0)=B_n ^{(\alpha)}$ are called the {\it{$n$-th Bernoulli number}} of order $\alpha$. By \eqref{1}, we easily get
\begin{equation*}
B_n ^{(\alpha)} (x)=\left(B^{(\alpha)} +x\right)^n=\sum_{l=0} ^n \binom{n}{l} B_l ^{(\alpha)} x^{n-l},
\end{equation*}
with the usual convention about replacing $(B^{(\alpha)})^n$ by $B_n ^{(\alpha)}$.

As is well known, the {\it{Euler polynomials}} of order $\alpha$ are also defined by the generating function to be
\begin{equation}\label{2}
\left(\frac{2}{e^t+1}\right)^{\alpha} e^{xt}=\sum_{n=0} ^{\infty}E_n ^{(\alpha)}(x)\frac{t^n}{n!},{\text{ (see [11,12,13])}}.
\end{equation}
In the special case, $x=0$, $E_n ^{(\alpha)} (0)=E_n ^{(\alpha)}$ are called the {\it{$n$-th Euler numbers}} of order $\alpha$. From \eqref{2}, we note that
\begin{equation*}
E_n ^{(\alpha)} (x)=\sum_{l=0} ^n \binom{n}{l} E_n ^{(\alpha)} x^{n-l},{\text{ (see [1-18])}}.
\end{equation*}
For $\lambda \in {\mathbb{C}}$ with $\lambda \neq 1$, the {\it{Frobenius-Euler polynomials}} of order $\alpha$ are defined by the generating function to be
\begin{equation}\label{3}
\left(\frac{1-\lambda}{e^t-\lambda}\right)^{\alpha} e^{xt}=\sum_{n=0} ^{\infty}H_n ^{(\alpha)}(x|\lambda)\frac{t^n}{n!},{\text{ (see [1,9,10,16])}}.
\end{equation}
In the special case, $x=0$, $H_n ^{(\alpha)} (0|\lambda)=H_n ^{(\alpha)} (\lambda)$ are called {\it{$n$-th Frobenius-Euler numbers}} of order $\alpha$. By \eqref{3}, we get
\begin{equation*}
H_n ^{(\alpha)} (x|\lambda)=\sum_{l=0} ^n \binom{n}{l}H_l(\lambda)x^{n-l},{\text{ (see [1,9,10,16])}}.
\end{equation*}
Let ${\mathbb{P}}$ be the algebra of polynomials in the variable $x$ over ${\mathbb{C}}$ and let ${\mathbb{P}}^{*}$ be the vector space of all linear functionals on ${\mathbb{P}}$. The action of the linear functional $L$ on a polynomial $p(x)$ is denoted by $\left<L|p(x)\right>$. We recall that the vector space on ${\mathbb{P}}^{*}$ are defined by $\left< L+M|p(x)\right>=\left<L|p(x)\right>+\left<M|p(x)\right>$, $\left< cL|p(x)\right>=c\left<L|p(x)\right>$, where $c$ is a complex constant.

Let
\begin{equation}\label{4}
{\mathcal{F}}=\left\{ \left.f(t)=\sum_{k=0} ^{\infty} \frac{a_k}{k!} t^k~\right|~ a_k \in {\mathbb{C}} \right\}.
\end{equation}
For $f(t)\in {\mathcal{F}}$, we define a linear functional on ${\mathbb{P}}$ by setting
\begin{equation}\label{5}
\left<f(t)|x^n \right>=a_n,{\text{ for all }}n\geq 0,{\text{ (see [13,15])}}.
\end{equation}
From \eqref{4} and \eqref{5}, we note that
\begin{equation}\label{6}
\left<t^k | x^n \right>=n! \delta_{n,k},~(n,k \geq 0), {\text{ (see [13,15])}},
\end{equation}
where $\delta_{n,k}$ is the Kronecker symbol.

For $f_L(t)=\sum_{k=0} ^{\infty} \frac{\left<L|x^k\right>}{k!}t^k$, we have $\left<f_L(t)|x^n\right>=\left<L|x^n\right>$. So the map $L \mapsto f_L(t)$ is a vector space isomorphism from ${\mathbb{P}}^{*}$ onto ${\mathcal{F}}$. Henceforth, ${\mathcal{F}}$ is thought of as both a formal power series and a linear functional. We shall call ${\mathcal{F}}$ the {\it{umbral algebra}}. The umbral calculus is the study of umbral algebra.

The order $o(f(t))$ of the non-zero power series $f(t)$ is the smallest integer $k$ for which the coefficient of $t^k$ does not vanish (see \cite{13}). If $o(f(t))=1$, then $f(t)$ is called a {\it{delta series}}. If $o(f(t))=0$, then $f(t)$ is called an {\it{invertible series}} (see \cite{15}). Let $o(f(t))=1$ and $o(g(t))=0$. Then there exists a unique sequence $S_n(x)$ of polynomials such that $\left<g(t)f(t)^k|S_n(x)\right>=n!\delta_{n,k}$ $(n,k \geq 0)$. The sequence $S_n(x)$ is called {\it{Sheffer sequence}} for $(g(t),f(t))$, which is denoted by $S_n(x)\sim (g(t),f(t))$ (see \cite{13,15}). Let $f(t)\in{\mathcal{F}}$ and $p(x)\in{\mathbb{P}}$. Then, by \eqref{6}, we easily see that $\left. \left<e^{yt}\right| p(x)\right>=p(y)$, $\left. \left<f(t)g(t)\right| p(x)\right>=\left. \left<f(t)\right| g(t)p(x)\right>=\left. \left<g(t)\right| f(t)p(x)\right>$.

For $f(t)\in{\mathcal{F}}$ and $p(x)\in{\mathbb{P}}$, we have
\begin{equation}\label{7}
f(t)=\sum_{k=0} ^{\infty} \frac{\left<f(t)|x^k\right>}{k!} t^k,~p(x)=\sum_{k=0} ^{\infty} \frac{\left<t^k|p(x)\right>}{k!}x^k{\text{ (see [13,15])}}.
\end{equation}
By \eqref{7}, we easily get
\begin{equation}\label{8}
p^{(k)}(0)=\left<t^k|p(x)\right>,~\left<1\left|p^{(k)}(x)\right.\right>=p^{(k)}(0).
\end{equation}
From \eqref{8}, we have
\begin{equation}\label{9}
t^kp(x)=p^{(k)}(x)=\frac{d^kp(x)}{dx^k},~(k \geq 0), {\text{ (see [13,5])}}.
\end{equation}

For $S_n(x) \sim \left(g(t),f(t)\right)$, the generating function of Sheffer sequence $S_n(x)$ is given by
\begin{equation}\label{10}
\frac{1}{g({\bar{f}}(t))}e^{y{\bar{f}}(t)}=\sum_{k=0} ^{\infty} \frac{S_k(y)}{k!}t^k,{\text{ for all }}y \in {\mathbb{C}},
\end{equation}
where ${\bar{f}}(t)$ is the compositional inverse of $f(t)$ (see \cite{15}).
Let us assume that
\begin{equation}\label{11}
S_n(x)\sim (1,f(t)),~t_n(x)\sim(1,g(t)).
\end{equation}
Then, we note that
\begin{equation}\label{12}
S_n(x)=x\left(\frac{g(t)}{f(t)}\right)^n x^{-1} t_n(x),{\text{ (see [13,15])}}.
\end{equation}
By \eqref{6}, we easily see that $x^n\sim(1,t)$.

In this paper, we give some interesting identities of higher-order Bernoulli, Frobenius-Euler and Euler polynomials involving multiple power and alternating sums which are derived from umbral calculus. By using our methods of this paper, we can obtain many interesting identities of special polynomials.

\section{Higher-order Bernoulli, Frobenius-Euler and Euler polynomials}

Let $S_n(x)\sim(g(t),f(t))$. Then we see that
\begin{equation}\label{13}
g(t)S_n(x) \sim (1,f(t)).
\end{equation}
From \eqref{12}, \eqref{13} and $x^n\sim(1,t)$, we note that
\begin{equation}\label{14}
S_n(x)=\frac{1}{g(t)} x \left(\frac{t}{f(t)}\right)x^{n-1}.
\end{equation}
The equation \eqref{14} is important to derive our results in this paper. From \eqref{1}, \eqref{2}, \eqref{3} and \eqref{10}, we can derive the following lemma:
\begin{lemma}\label{lem1}
For $n \geq 0$, $m \in {\mathbb{N}}$, we have
\begin{equation*}
\begin{array}{cc}
m^n B_n ^{(\alpha)} \left(\frac{x}{m}\right)\sim\left(\left(\frac{e^{mt}-1}{mt}\right)^{\alpha},t\right),~\frac{m^n}{m^{\alpha}}B_n ^{(\alpha)}\left(\frac{x}{m}\right)\sim\left(\left(\frac{e^{mt}-1}{t}\right)^{\alpha},t\right), \\
m^n E_n ^{(\alpha)} \left(\frac{x}{m}\right)\sim\left(\left(\frac{e^{mt}+1}{2}\right)^{\alpha},t\right),~m^nH_n ^{(\alpha)}\left.\left(\frac{x}{m}\right|\lambda\right)\sim\left(\left(\frac{e^{mt}-\lambda}{1-\lambda}\right)^{\alpha},t\right).
\end{array}
\end{equation*}
\end{lemma}
Let us consider the following Sheffer sequences:
\begin{equation}\label{15}
S_n(x)\sim\left(1,\frac{t^2}{e^t-1}\right),~t_n(x)\sim\left(1,\frac{t^2}{e^{mt}-1}\right).
\end{equation}
From \eqref{14}, we have
\begin{equation}\label{16}
\begin{split}
S_n(x)=&x\left(\frac{e^t-1}{t}\right)^nx^{n-1}=x\sum_{l=0} ^{\infty}\frac{n!}{(l+n)!}S_2(l+n,n)t^lx^{n-1}\\
=&x\sum_{l=0} ^{n-1} \frac{\binom{n-1}{l}}{\binom{l+n}{n}}S_2(l+n,n)x^{n-1-l}\\
=&x\sum_{r=0} ^{n-1} \frac{\binom{n-1}{r}}{\binom{2n-1-r}{n}}S_2(2n-1-r,n)x^r,
\end{split}
\end{equation}
and
\begin{equation}\label{17}
\begin{split}
t_n(x)=&x\left(\frac{e^{mt}-1}{t}\right)^nx^{n-1}=x\sum_{l=0} ^{\infty}\frac{n!}{(l+n)!}S_2(l+n,n)m^{l+n}t^lx^{n-1}\\
=&x\sum_{l=0} ^{n-1} \frac{n!}{(l+n)!}S_2(l+n,n)m^{n+l}(n-1)_lx^{n-1-l} \\
=&x\sum_{l=0} ^{n-1} \frac{\binom{n-1}{l}}{\binom{l+n}{n}}S_2(l+n,n)m^{n+l}x^{n-1-l}\\
=&x\sum_{r=0} ^{n-1} \frac{\binom{n-1}{r}}{\binom{2n-1-r}{n}}S_2(2n-1-r,n)m^{2n-1-r}x^r,
\end{split}
\end{equation}
where $S_2(n,k)$ is the Stirling number of the second kind.

For $n \geq 1$, by \eqref{12} and \eqref{15}, we get
\begin{equation}\label{18}
\begin{split}
t_n(x)=&x\left(\frac{e^{mt}-1}{e^t-1}\right)^nx^{-1}S_n(x)=x\left(e^{-t}\sum_{l=1} ^m e^{lt}\right)^nx^{-1}S_n(x) \\
=&xe^{-nt}\sum_{\begin{array}{cc}
0 \leq v_1,\ldots,v_m \leq n \\
v_1+\cdots+v_m=n
\end{array}}\binom{n}{v_1,\ldots,v_m}e^{(v_1+2v_2+\cdots+mv_m)t}x^{-1}S_n(x) \\
=&x\left.\sum_{s=0} ^{\infty}\right\{\sum_{k=0} ^s \binom{s}{k}(-n)^{s-k}\sum_{\begin{array}{cc}
0 \leq v_1,\ldots,v_m \leq n \\
v_1+\cdots+v_m=n
\end{array}}\binom{n}{v_1,\ldots,v_m} \\
& \times \left.(v_1+2v_2+\cdots+mv_m)^k \frac{t^s}{s!}\right\}x^{-1}S_n(x).
\end{split}
\end{equation}
Let us define multiple power sum $S_k ^{(n)} (m)$ as follows:
\begin{equation}\label{19}
S_k ^{(n)} (m)=\sum_{\begin{array}{cc}
0 \leq v_1,\ldots,v_m \leq n \\
v_1+\cdots+v_m=n
\end{array}}\binom{n}{v_1,\ldots,v_m}(v_1+2v_2+\cdots+mv_m)^k.
\end{equation}
By \eqref{17}, \eqref{18} and \eqref{19}, we get
\begin{equation}\label{20}
\begin{split}
t_n(x)=&x\sum_{s=0} ^{\infty} \sum_{k=0} ^s \binom{s}{k}(-n)^{s-k}S_k ^{(n)} (m) \frac{t^s}{s!}x^{-1}S_n(x) \\
=&x\sum_{r=0} ^{n-1} \sum_{s=0} ^r \sum_{k=0} ^s \frac{\binom{r}{s}\binom{n-1}{r}\binom{s}{k}}{\binom{2n-1-r}{n}}(-n)^{s-k}S_2(2n-1-r,n)S_k ^{(n)} (m)x^{r-s}.
\end{split}
\end{equation}
From \eqref{12} and \eqref{15}, we can also derive
\begin{equation}\label{21}
\begin{split}
t_n(x)=&x\left(\frac {e^{mt}-1}{e^t-1}\right)^nx^{-1}S_n(x)=x\left(\frac{e^{mt}-1}{t}\right)^n\left(\frac{t}{e^t-1}\right)^nx^{-1}S_n(x) \\
=&x\left(\frac {e^{mt}-1}{t}\right)^n\sum_{r=0} ^{n-1} \frac{\binom{n-1}{r}}{\binom{2n-1-r}{n}}S_2(2n-1-r,n)B_r ^{(n)}(x)\\
=&x\sum_{r=0} ^{n-1} \sum_{s=0} ^r \frac{\binom{r}{s}\binom{n-1}{r}}{\binom{s+n}{n}\binom{2n-1-r}{n}}S_2(s+n,n)S_2(2n-1-r,n)m^{n+s}B_{r-s} ^{(n)}(x).
\end{split}
\end{equation}
Therefore, by \eqref{20} and \eqref{21}, we obtain the following theorem.
\begin{theorem}\label{thm3}
For $n \geq 1$, we have
\begin{equation*}
\begin{split}
&\sum_{r=0} ^{n-1} \sum_{s=0} ^r \sum_{k=0} ^s \frac{\binom{r}{s}\binom{n-1}{r}\binom{s}{k}}{\binom{2n-1-r}{n}}(-n)^{s-k}S_2(2n-1-r,n)S_k ^{(n)} (m)x^{r-s} \\
=&\sum_{r=0} ^{n-1} \sum_{s=0} ^r \frac{\binom{r}{s}\binom{n-1}{r}}{\binom{s+n}{n}\binom{2n-1-r}{n}}S_2(s+n,n)S_2(2n-1-r,n)m^{n+s}B_{r-s} ^{(n)} (x).
\end{split}
\end{equation*}
\end{theorem}
Let us consider the following Sheffer sequences:
\begin{equation}\label{22}
S_n(x)\sim\left(1,e^t-1\right),~t_n(x)\sim\left(1,e^{mt}-1\right),
\end{equation}
where $m \in {\mathbb{N}}$ and $n \geq 0$.

For $n \geq 1$, by \eqref{14}, we get
\begin{equation}\label{23}
S_n(x)=x\left(\frac{t}{e^t-1}\right)^n x^{n-1}=xB_{n-1} ^{(n)} (x),
\end{equation}
and
\begin{equation}\label{24}
t_n(x)=x\left(\frac{t}{e^{mt}-1}\right)^nx^{n-1}.
\end{equation}
By Lemma \ref{lem1} and \eqref{24}, we get
\begin{equation}\label{25}
t_n(x)=\frac{x}{m}B_{n-1} ^{(n)} \left(\frac{x}{m}\right).
\end{equation}
From \eqref{12} and \eqref{22}, we can derive
\begin{equation}\label{26}
\begin{split}
S_n(x)=&x\left(\frac{e^{mt}-1}{e^t-1}\right)^nx^{-1}t_n(x)=x\left(e^{-t}\sum_{l=0} ^m e^{lt}\right)^n x^{-1}t_n(x)\\
=&\frac{x}{m}\sum_{s=0} ^{\infty} \sum_{k=0} ^s \binom{s}{k} (-n)^{s-k} S_k ^{(n)}(m)\frac{t^s}{s!}B_{n-1} ^{(n)} \left(\frac{x}{m}\right)\\
=&x\sum_{s=0} ^{n-1} \sum_{k=0} ^s \binom{s}{k}\binom{n-1}{s}(-n)^{s-k}S_k ^{(n)} (m)B_{n-1-s} ^{(n)} \left(\frac{x}{m}\right)m^{-s-1}.
\end{split}
\end{equation}
\begin{theorem}\label{thm4}
For $n,m\geq 1$, we have
\begin{equation*}
B_{n-1} ^{(n)}(x)=\\ \sum_{r=0} ^{n-1} \sum_{s=0} ^r \sum_{k=0} ^s \frac{\binom{r}{s}\binom{n-1}{r}\binom{s}{k}}{\binom{2n-1-r}{n}}(-n)^{s-k}S_2(2n-1-r,n)S_k ^{(n)} (m)x^{r-s}.
\end{equation*}
\end{theorem}
Let us assume that
\begin{equation}\label{27}
S_n(x)\sim\left(1,t\left(\frac{e^t+1}{2}\right)\right),~t_n(x)\sim\left(1,\left(\frac{e^{mt}+1}{2}\right)t\right),
\end{equation}
where $n \geq 0$ and $m \in {\mathbb{N}}$ with $m \equiv1$ $({\rm{mod}}~2)$. By \eqref{14}, we get
\begin{equation}\label{28}
S_n(x)=x\left(\frac{2}{e^t+1}\right)^nx^{n-1}=xE_{n-1} ^{(n)} (x),
\end{equation}
and
\begin{equation}\label{29}
t_n(x)=x\left(\frac{2}{e^{mt}+1}\right)^nx^{n-1}.
\end{equation}
From Lemma \ref{lem1} and \eqref{29}, we note that
\begin{equation}\label{30}
t_n(x)=xm^{n-1}E_{n-1} ^{(n)} \left(\frac{x}{m}\right).
\end{equation}
By \eqref{12} and \eqref{27}, we see that
\begin{equation}\label{31}
\begin{split}
&S_n(x)\\
=&x\left(\frac{e^{mt}+1}{e^t+1}\right)^nx^{-1}t_n(x)=x\left(-e^{-t}\sum_{l=1} ^m (-e^t)^l \right)^nx^{-1}t_n(x)\\
=&(-1)^nxe^{-nt}\sum_{\begin{array}{cc}
0\leq v_1,\ldots,v_m \leq n\\
v_1+\cdots+v_m=n
\end{array}}
\binom{n}{v_1,\ldots,v_m} \\
& \times (-1)^{v_1+2v_2+\cdots+mv_m}e^{(v_1+2v_2+\cdots+mv_m)t}x^{-1}t_n(x)\\
=&(-1)^nxe^{-nt}\left.\sum_{k=0} ^{\infty}\right(\sum_{\begin{array}{cc}
0\leq v_1,\ldots,v_m \leq n\\
v_1+\cdots+v_m=n
\end{array}}
\binom{n}{v_1,\ldots,v_m} \\
& \times (-1)^{v_1+2v_2+\cdots+mv_m}\left(v_1+2v_2+\cdots+mv_m\right)^k\left)\frac{t^k}{k!}\right.x^{-1}t_n(x).
\end{split}
\end{equation}
Let us define multiple alternating power sums $T_k ^{(n)} (m)$ as follows:
\begin{equation}\label{32}
\begin{split}
T_k ^{(n)} (m)=&\sum_{\begin{array}{cc}
0\leq v_1,\ldots,v_m \leq n\\
v_1+\cdots+v_m=n
\end{array}}
\binom{n}{v_1,\ldots,v_m}  \\
& \times (-1)^{v_1+2v_2+\cdots+mv_m}\left(v_1+2v_2+\cdots+mv_m\right)^k.
\end{split}
\end{equation}
By \eqref{31} and \eqref{32}, we get
\begin{equation}\label{33}
\begin{split}
S_n(x)=&(-1)^nx\sum_{s=0} ^{\infty}\sum_{k=0} ^s \binom{s}{k}(-n)^{s-k}T_k ^{(n)} (m)\frac{t^s}{s!}m^{n-1}E_{n-1} ^{(n)}\left(\frac{x}{m}\right)\\
=&(-1)^nx\sum_{s=0} ^{n-1}\sum_{k=0} ^s \binom{s}{k}\binom{n-1}{s}(-n)^{s-k}T_k ^{(n)}(m)m^{n-s-1}E_{n-1-s} ^{(n)} \left(\frac{x}{m}\right).
\end{split}
\end{equation}
Therefore, by \eqref{28} and \eqref{33}, we obtain the following theorem.
\begin{theorem}\label{thm5}
For $n,m \geq 1$ with $m \equiv1~({\rm{mod}}~2)$, we have
\begin{equation*}
E_{n-1} ^{(n)} (x)=(-1)^n\sum_{s=0} ^{n-1}\sum_{k=0} ^s \binom{s}{k}\binom{n-1}{s}(-n)^{s-k}T_k ^{(n)}(m)m^{n-s-1}E_{n-1-s} ^{(n)} \left(\frac{x}{m}\right).
\end{equation*}
\end{theorem}
Let us consider the following Sheffer sequences:
\begin{equation}\label{34}
S_n(x)\sim\left(1,t\left(\frac{e^t-\lambda}{1-\lambda}\right)\right),~t_n(x)\sim\left(1,t\left(\frac{e^{mt}-\lambda^m}{1-\lambda^m}\right)\right),
\end{equation}
where $m \in {\mathbb{N}}$ and $\lambda \in {\mathbb{C}}$ with $\lambda \neq 0$, $\lambda^m \neq 1$. By \eqref{14}, we get
\begin{equation}\label{35}
S_n(x)=x\left(\frac{1-\lambda}{e^t-\lambda}\right)^nx^{n-1}=xH_{n-1} ^{(n)} (x|\lambda),
\end{equation}
and
\begin{equation}\label{36}
t_n(x)=x\left(\frac{1-\lambda^m}{e^{mt}-\lambda^m} \right)^{n} x^{n-1}.
\end{equation}
From Lemma \ref{lem1} and \eqref{36}, we have
\begin{equation}\label{37}
t_n(x)=m^{n-1}xH_{n-1} ^{(n)} \left.\left(\frac{x}{m}\right|\lambda^m\right).
\end{equation}
By \eqref{12} and \eqref{34}, we get
\begin{equation}\label{38}
\begin{split}
&S_n(x)\\
=&x\left(\frac{1-\lambda}{1-\lambda^m}\right)^n\left(\frac{e^{mt}-\lambda^m}{e^t-\lambda}\right)^nx^{-1}t_n(x)\\
=&x\left(\frac{1-\lambda}{1-\lambda^m}\right)^n\lambda^{mn-n}\left(\frac{1-\left(\frac{e^t}{\lambda}\right)^m}{1-\frac{e^t}{\lambda}}\right)^nx^{-1}t_n(x)\\
=&x\left(\frac{1-\lambda}{1-\lambda^m}\right)^n\lambda^{mn-n}\left(\frac{\lambda}{e^t}\sum_{l=1} ^m \left(\frac{e^t}{\lambda}\right)^l\right)^nx^{-1}t_n(x)\\
=&\left(\frac{1-\lambda}{1-\lambda^m}\right)^n\lambda^{mn}xe^{-nt}\sum_{\begin{array}{cc}
0\leq v_1,\ldots,v_m\leq n\\
v_1+\cdots+v_m=n
\end{array}}\binom{n}{v_1,\ldots,v_m} \\
&\times \lambda^{-(v_1+2v_2+\cdots+mv_m)}e^{(v_1+2v_2+\cdots+mv_m)t}x^{-1}t_n(x)\\
=&\left(\frac{1-\lambda}{1-\lambda^m}\right)^n\lambda^{mn}x\sum_{s=0} ^{\infty}\sum_{k=0} ^s \binom{s}{k}(-n)^{s-k}\sum_{\begin{array}{cc}
0\leq v_1,\ldots,v_m\leq n\\
v_1+\cdots+v_m=n
\end{array}}\binom{n}{v_1,\ldots,v_m} \\
&\times \lambda^{-(v_1+2v_2+\cdots+mv_m)}(v_1+2v_2+\cdots+mv_m)^k\frac{t^s}{s!}m^{n-1}H_{n-1} ^{(n)} \left.\left(\frac{x}{m}\right|\lambda^m\right).
\end{split}
\end{equation}
Let us define $\lambda$-analogue of multiple power sums $S_k ^{(n)}(m|\lambda)$ as follows:
\begin{equation}\label{39}
\begin{split}
S_k ^{(n)} (m|\lambda)=&\sum_{\begin{array}{cc}
0\leq v_1,\ldots,v_m\leq n\\
v_1+\cdots+v_m=n
\end{array}}\binom{n}{v_1,\ldots,v_m}  \\
& \times \lambda^{-(v_1+2v_2+\cdots+mv_m)}(v_1+2v_2+\cdots+mv_m)^k.
\end{split}
\end{equation}
From \eqref{38} and \eqref{39}, we have
\begin{equation}\label{40}
\begin{split}
\aligned
S_n(x)=&\left(\frac{1-\lambda}{1-\lambda^m}\right)^n\lambda^{mn}x\sum_{s=0} ^{n-1}\sum_{k=0} ^s\binom{s}{k}(-n)^{s-k}S_k ^{(n)} (m|\lambda)\frac{t^s}{s!}m^{n-1}H_{n-1} ^{(n)}\left.\left(\frac{x}{m}\right|\lambda^m\right) \\
=&\left(\frac{1-\lambda}{1-\lambda^m}\right)^n\lambda^{mn}x\sum_{s=0} ^{n-1}\sum_{k=0} ^s\binom{s}{k}(-n)^{s-k}S_k ^{(n)} (m|\lambda)
\\ &\times \binom{n-1}{s}m^{n-1-s}H_{n-1-s} ^{(n)} \left.\left(\frac{x}{m}\right|\lambda^m\right).
\endaligned
\end{split}
\end{equation}
Therefore, by \eqref{35} and \eqref{40}, we obtain the following theorem.
\begin{theorem}
For $m,n \geq 1$, $\lambda \in {\mathbb{C}}$ with $\lambda \neq 0$, $\lambda^m \neq 1$ , we have
\begin{equation*}
\aligned
 H_{n-1} ^{(n)}(x|\lambda)=& \left(\frac{1-\lambda}{1-\lambda^m}\right)^n\lambda^{mn}\sum_{s=0} ^{n-1}\sum_{k=0} ^s\binom{s}{k}\binom{n-1}{s}(-n)^{s-k}m^{n-1-s}
\\ & \times S_k ^{(n)} (m|\lambda)  H_{n-1-s} ^{(n)} \left.\left(\frac{x}{m}\right|\lambda^m\right).
\endaligned
\end{equation*}

\end{theorem}

\end{document}